\documentclass[a4paper,numbers=noendperiod]{scrartcl}

\usepackage[english]{babel}
\usepackage[utf8]{inputenc}
\usepackage[T1]{fontenc}
\usepackage{lmodern}

\usepackage{multirow}
\usepackage{rotating}
\usepackage{booktabs}
\usepackage{tabularx}
\newcolumntype{C}{>{\centering\arraybackslash}X}
\newcolumntype{R}{>{\raggedleft\arraybackslash}X}
\newcolumntype{L}{>{\raggedright\arraybackslash}X}

\usepackage{algorithm}
\usepackage{algorithmic}
\renewcommand{\algorithmicrequire}{\textbf{Input:}}

\usepackage[usenames,dvipsnames,table]{xcolor}
\definecolor{rgbcyan}{HTML}{00ADEF}
\definecolor{rgbmagenta}{HTML}{EC008C}

\usepackage{tikz}
\usetikzlibrary{mindmap,backgrounds}
\usetikzlibrary{shapes,arrows}
\usetikzlibrary{matrix}
\usetikzlibrary{positioning}
\usetikzlibrary{calc}
\tikzstyle{line} = [draw, thick]
\usepackage[font=footnotesize,skip=4pt]{subcaption}

\graphicspath{{../}}

\usepackage{pifont}
\usepackage{amsmath}
\usepackage{amssymb}
\usepackage{amsthm}

\newtheorem{definition}{Definition}

\renewcommand{\vec}[1]{\boldsymbol{#1}}

\usepackage{hyperref}

\usepackage{siunitx}

\title{The True Destination of EGO is Multi-local Optimization}

\author{Simon Wessing\textsuperscript{1}, Mike Preuss\textsuperscript{2}\\\textsuperscript{1}Computer Science Department\\Technische Universität Dortmund, Germany\\
\textsuperscript{2}European Research Center for Information Systems\\WWU M{\"u}nster, Germany\\
\url{simon.wessing@tu-dortmund.de}, \url{mike.preuss@uni-muenster.de}}

\date{}
\begin{document}
 
\maketitle

\begin{abstract}
Efficient global optimization is a popular algorithm for the optimization of expensive multimodal black-box functions. 
One important reason for its popularity is its theoretical foundation of global convergence. 
However, as the budgets in expensive optimization are very small, the asymptotic properties only play a minor role and the algorithm sometimes comes off badly in experimental comparisons.
Many alternative variants have therefore been proposed over the years.
In this work, we show experimentally that the algorithm instead has its strength in a setting where multiple optima are to be identified.
\end{abstract}

\section{Introduction}

Efficient global optimization (EGO) is a popular algorithm for the optimization of expensive multimodal black-box functions. 
At its core is a Kriging metamodel, whose predictions are used to formulate a so-called infill criterion. 
This criterion usually is a compromise between two goals: a) to detect especially good solutions, and b) to improve the model itself in order to enable better predictions.
The established model can be employed to cheaply search for potential new points because it is much faster than the original function, often by a factor of $1000$ or more.
By optimizing with regard to the infill criterion, a new point can be determined for sampling the expensive function. 
It is clear that this works reasonably well only for functions that can be predicted from a sparse sample, namely relatively low dimensional and generally rather smooth ones. 
With more and more samples coming in, the model improves, so that one gets a better and better overall impression of the original function. 
However, the model fit also gets more and more expensive, due to necessary matrix inversions.
Thus, the sample size is limited, usually to around $1000$ points.

The infill criterion of choice for EGO is the so-called expected improvement (EI), which incorporates both objectives mentioned above.
While EGO is conceptually elegant and its convergence rate to the global optimum can be analyzed mathematically~\cite{Bull2011}, there are multiple experimental results that indicate a preference for a more greedy infill criterion in expensive global optimization.
For example, S\'obester et al.~\cite{Sobester2005} propose a weighted expected improvement, which can be adjusted to search more locally or globally, depending on the weights.
Some more evidence has surfaced in research on ``multipoint'' infill criteria for parallelizing function evaluations.
Bischl et al.~\cite{Bischl2014} discovered that just using the model prediction was the most successful infill criterion in a comparison with expected improvement and several multipoint infill criteria for a budget of $45n$ objective function evaluations, where $n$ is the number of decision variables of the problem.
Ginsbourger et al.~\cite{Ginsbourger2010} showed that an explorative variant of their \emph{constant liar} criterion was less competitive than the more exploitative one, when filling in four to ten points on the Branin function. 
Also Ursem~\cite{Ursem2014} developed an \emph{investment portfolio improvement} function to propose three solutions per iteration, ranging from high exploitation to high exploration.

Our position is, that while the rather global search strategy of expected improvement may prevent a highly accurate approximation of the global optimum with small budgets, it represents a virtue for the task of finding multiple optima. 
This task is also known under the names of multi-global or multi-local optimization, depending on if only global or all optima are sought.
EGO has to store the sampled points and function values anyway, to build the model, and all we have to do is to add a basin identification heuristic, to decide which points correspond to distinct attraction basins, and simply select the best one of each basin. 
As we want to avoid the responsibility of deciding if an optimum is global or local in our algorithm, we focus on multi-local optimization here.

Employing EGO as multi-local optimization algorithm may seem counterintuitive at first, but we claim that in an expensive optimization scenario, where we can afford only several hundreds of function evaluations, this makes a lot of sense. 
Many competing multi-local optimization algorithms rely on (multiple) local searches which turn out to be too expensive in this scenario.
The focus of this paper is to experimentally analyse how well EGO is suited for expensive multi-local optimization. 
As expensive optimization setting, we assume budgets of $500$ evaluations or less here. 
With respect to the limitation concerning the number of points that can be used as basis of a Kriging model, EGO matches the expensive optimization scenario very well. 
We therefore run EGO against one of the most effective and robust black-box optimization methods, the (Restart-) CMA-ES~\cite{Auger2005a}. 
The CMA-ES is a good reference as it is also part of modern multimodal optimization algorithms, e.\,g., it performs the local optimization step in the NEA2 method \cite{Preuss2015}. 
Feliot et al.~\cite{Feliot2017} also compared model-based approaches with local search algorithms in a constrained optimization setting and found that the competition is quite balanced. 
This already shows that the situation is not as clear-cut as it seems.

\begin{figure*}[t]
    \centering
    \includegraphics[scale=0.495]{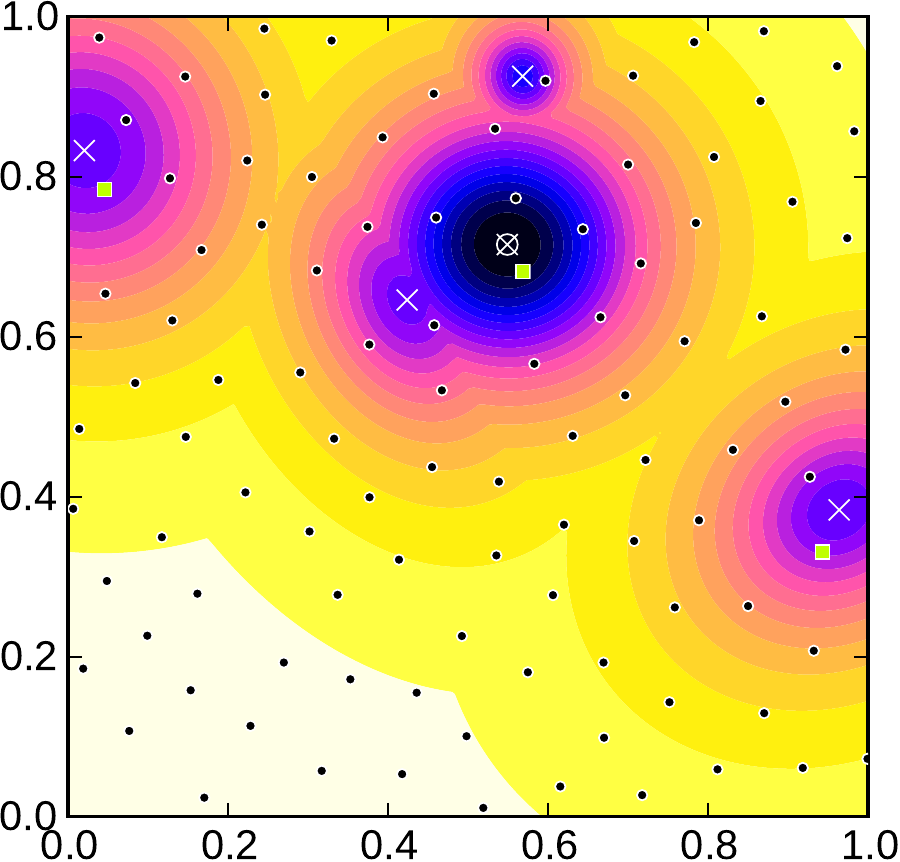}\hfill%
\includegraphics[scale=0.495]{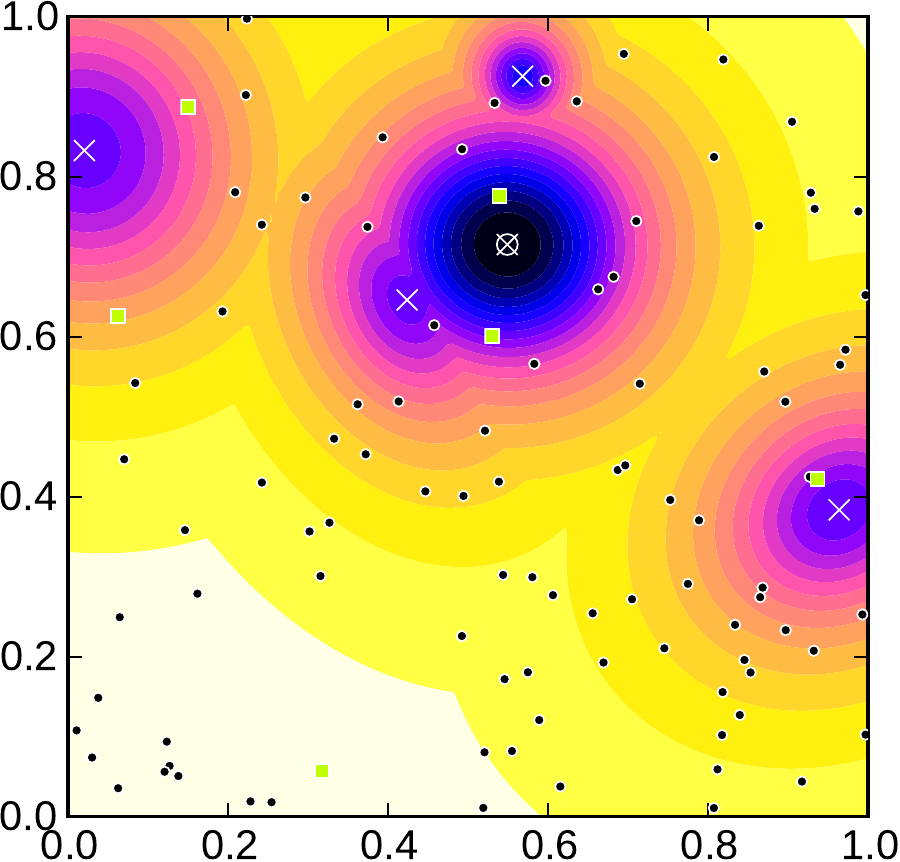}\hfill%
\includegraphics[scale=0.495]{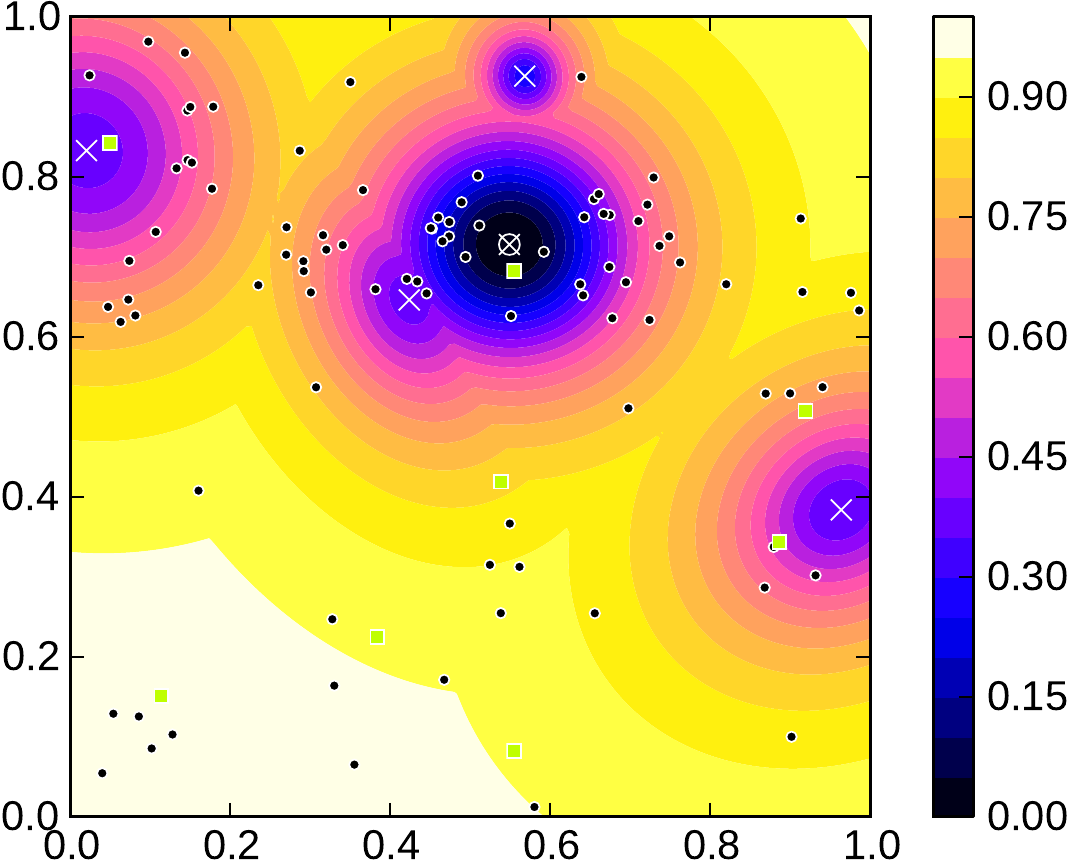}
    \caption{Examples of nearest-better clustering on a 2-D landscape, applied to $N = 100$ highly uniform, random uniform, and clustered points (from left to right). Selected points representing basin centers are marked with green squares.}
	\label{fig:nbc_failure_example}
\end{figure*}

Our comparison is carried out on a mix of typical global optimization benchmark functions and niching competition benchmark functions, and we look at the results from two perspectives: 
\begin{enumerate}
\item global optimization: we are only interested in locating the global optimum (one global optimizer in case it has several preimages), 
\item multi-local optimization: we want to detect as many local/global optima as possible, ideally all of them.
\end{enumerate}
If our reasoning from above is correct, EGO should perform better than CMA-ES and related methods under a multi-local optimization perspective, but worse if seen from a global optimization perspective. This would mean that in the expensive black-box setting, EGO is a very suitable multi-local optimization method. In order to compare it to other such methods, we need to add a basin identification method, because in the answer set of the algorithm, we only want points that are approximations of existing optima, not the full archive. According to~\cite{Wessing2016a}, we have two methods at hand: topographical selection (TS) by~\cite{Toern1992} and nearest-better clustering (NBC) as described in~\cite{PreussNiching2012}. 

Previous investigations showed that using NBC can be problematic because it produces too many clusters if the point set deviates from uniformity. 
An example for this effect can be found in Figure~\ref{fig:nbc_failure_example}, where the uniformity decreases from left to right.
If, e.\,g., the clustered points reside on a linear slope, the best point of a cluster is wrongly interpreted as representing an optimum where none exists (see rightmost subfigure, bottom left corner).
The reason for this behavior is that outliers tend to be selected simply because of their large nearest-neighbor distances. 
We therefore rely on TS in this work.

To our knowledge, EGO and related surrogate-model-based methods have never been investigated in this way.  
Also the whole area of expensive multi-local/multimodal optimization seems to have been sparsely visited.
One of the few works in this area is~\cite{Parmee1996}, however, they do not employ the term \emph{multimodal optimization} as it has been established only years later. 

The next section introduces necessary problem definitions, Sect.~\ref{sec:filter_methods} explains the employed methods, namely EGO and topographical selection. 
The remainder of the paper is mostly concerned with the comparison experiment (Sect.~\ref{sec:experiment}) and ends with the conclusions.

\section{Problem definition}
\label{sec:problem-definition}

In the following, we will assume to have a deterministic objective function $f: \mathcal{X} \to \mathbb{R}$, where $\mathcal{X} = [\vec{\ell}, \vec{u}] \subset \mathbb{R}^n$ is the \emph{search space} or \emph{region of interest} (ROI) and $n \in \mathbb{N}$ is the fixed number of decision variables.
The vectors $\vec{\ell} = (\ell_1, \dots, \ell_n)^\top$ and $\vec{u} = (u_1, \dots, u_n)^\top$ are called the lower and upper bounds of $\mathcal{X}$, respectively.
Let $N(\vec{y}) = \{\vec{x} \in \mathcal{X} \mid d(\vec{x}, \vec{y}) \leq \epsilon\}$ be the neighborhood of a point $\vec{y} \in \mathcal{X}$.
We say $f^* := f(\vec{x}^*)$ is a local minimum if $\exists \epsilon > 0: \nexists \vec{x} \in N(\vec{x}^*): f(\vec{x}) < f(\vec{x}^*)$. 
Technically, this implies that also plateaus are considered as local optima, although they are rather not intuitively, but at least this ensures that every position of a global optimum is also one of a local optimum.
A multimodal optimization problem was implicitly already defined by T\"orn and \v{Z}ilinskas~\cite[pp.~2--3]{Toern1989}. 
In~\cite[p.~6]{Wessing2015}, the definition was formulated explicitly as follows:

\begin{definition}[Multimodal minimization problem]
Let there be $\nu$ local minima $f_1^*, \dots, f_\nu^*$ of $f$ in $\mathcal{X}$.
If the ordering of these optima is $f_{(1)}^* < \dots < f_{(l)}^* < h < \dots < f_{(\nu)}^*$, a multimodal minimization problem is given as the task to approximate the set $\bigcup_{i=1}^l X_{(i)}^*$, where $X_{(i)}^* = \{\vec{x} \in \mathcal{X} \mid f(\vec{x}) = f_{(i)}^*\}$.
\end{definition}

The variable $h$ in this definition is simply a threshold to potentially exclude some of the worse optima. 
If $h = \infty$, we will be interested in all local optima of a problem. 
If $f_{(1)}^* < h < f_{(2)}^*$, we are only interested in approximating all the global optima.
Let $P$ be the obtained approximation set.
Additional constraints may be applied to this set, to obtain more specific problem definitions. 
For example, the cardinality of $P$ could be restricted by requiring $|P| \leq k$.
If $k = 1$, we have the conventional global optimization problem, where typically only one solution is sought.
Another issue are  diversity requirements, which could be formulated by demanding $\forall \vec{x}, \vec{y} \in P, \vec{x} \neq \vec{y}: d(\vec{x}, \vec{y}) > \epsilon$, i.\,e., the distance between any two solutions may not be smaller than some threshold $\epsilon$.
Alternatively, also more sophisticated diversity measures on the set $P$ may be calculated, maybe leading to a multiobjective formulation of the problem~\cite{Preuss2013_noeditors}.
However, we will stick to the basic task of finding all optima ($h = \infty$) here.

\section{Methods}
\label{sec:filter_methods}

\begin{algorithm}[tb]
\begin{algorithmic}[1]
\STATE generate an initial design $D \subset \mathcal{X}$\;
\STATE $\vec{y} \gets f(D)$\;
\WHILE{total evaluation budget is not exceeded}
      \STATE fit surrogate on $D$ and obtain $\hat{f}$, $\hat{s}$\;
      \STATE get new design point $\vec{x}'$ by optimizing the infill criterion based on $\hat{f}$, $\hat{s}$\;
      \STATE $y' \gets f(\vec{x}')$\COMMENT{evaluate new point}
      \STATE $D \leftarrow (D, \vec{x}')$\COMMENT{update design}
	  \STATE $\vec{y} \leftarrow (\vec{y}, y')$\COMMENT{update responses~}
    \ENDWHILE
    \RETURN $\hat{y}^* = \min(\vec{y})$ and the associated $\hat{\vec{x}}^*$
\end{algorithmic}
    \caption{Sequential model-based optimization}
    \label{alg:mbo}
\end{algorithm}

Algorithm~\ref{alg:mbo} illustrates the general sequential model-based optimization (MBO) framework, of which EGO is an instantiation. 
The main idea in model-based optimization is to approximate the expensive function $f(\vec{x})$ in
every iteration by a regression model, which is much cheaper to evaluate.
This is also called a meta-model or surrogate.
We are using a Kriging model, which not only provides a direct estimation $\hat{f}(\vec{x})$
of the true function value $f(\vec{x})$
but also an estimation of the prediction standard error $\hat{s}(\vec{x})$, also called a local uncertainty measure.
Our Kriging implementation follows the ``empirical Bayes'' approach with a correlation kernel and a maximum likelihood estimation of its parameters~\cite{Jones2001}.

The whole MBO concept has roots in response surface methodology, which was originally applied to physical experiments (with a human in the loop) \cite[pp.~7--11]{Castillo2007}.
It starts by exploring the parameter space with an initial design, often constructed in a
space-filling fashion. The main sequential loop can be divided into two alternating stages:
\begin{enumerate}
\item Fit a response surface to training data (including estimation of the model's parameters).
\item Use the surface to compute new search points under the assumption that the parameters are correct.
\end{enumerate}
In~\cite{Jones1998}, the now standard expected improvement criterion was proposed.
It is defined as
\begin{align*}
\text{EI}(\vec{x}) &= \mathbb{E}[\max \{0, \hat{y}^* - \hat{f}\}] = \left(\hat{y}^* - \hat{f}(\vec{x})\right) ~\Phi\left(\dfrac{\hat{y}^* - \hat{f}(\vec{x})}{\hat{s}(\vec{x})}\right) + \hat{s}(\vec{x}) ~\phi\left(\dfrac{\hat{y}^* - \hat{f}(\vec{x})}{\hat{s}(\vec{x})}\right) \,,
\end{align*}
where $\phi$ and $\Phi$ are the density and cumulative distribution function
of the standard normal distribution, respectively.
Hence, the sought point is $\vec{x}^* = \arg\max_{\vec{x} \in \mathcal{X}} \text{EI}(\vec{x})$.

\begin{algorithm}[t]
\renewcommand{\algorithmicrequire}{\textbf{Input:}}
\begin{algorithmic}[1]
\REQUIRE points $\mathcal{P} = \{\vec{x}_1, \dots, \vec{x}_N\}$, number $k$ of nearest neighbors
\ENSURE nodes of the topograph with no outgoing edges
\STATE create a directed graph $G = (V, E)$ with $V = \{v_1, \dots, v_N\}$ and $E = \emptyset$
\FORALL{$i \in \{1, \dots, N\}$}
  \STATE $J \gets$ indices of the $k$ nearest neighbors of $\vec{x}_i$ 
  \FORALL{$j \in J$}
  	\IF{$f(\vec{x}_j) < f(\vec{x}_i)$}
      \STATE $E \gets E \cup \{(v_i, v_j)\}$ \COMMENT{add edge to graph}
    \ELSIF{$f(\vec{x}_i) < f(\vec{x}_j)$}
      \STATE $E \gets E \cup \{(v_j, v_i)\}$ \COMMENT{add edge to graph}
  	\ENDIF
  \ENDFOR
\ENDFOR
\RETURN{$\{v \in V \mid \operatorname{deg}^+(v) = 0\}$} \COMMENT{select nodes with no outgoing edges}
\end{algorithmic}
\caption{Topographical selection}
\label{alg:ts}
\end{algorithm}

Topographical selection is provided as pseudo-code in algorithm~\ref{alg:ts}. It has similarities to nearest-better  clustering because the basic idea is to compare points regarding their quality to their closest neighbors. While NBC argues with relative distances, TS relies on fixed-size $k$-neighborhoods, such that $k$ is a parameter of the algorithm. 
We start with an empty graph and for every point, and detect the $k$ nearest neighbors. 
For each neighbor, we add an edge that points from the worse to the better point. 
After finishing the loop, we return all points that have only incoming edges as basin representatives. 

As performance measure for an approximation set $P$ in the case of global optimization, we use the deviation from the global optimum $f_\Delta = \hat{y}^* - f_1^*$.
In the multi-local case the number of found optima 
$
o = |\{\vec{x}^* \in X^* \mid d_{\mathrm{nn}}(\vec{x}^*, P) \leq r\}|
$
divided by the total number of optima $|X^*| = \nu$ as peak ratio $\operatorname{PR}(P) = o/\nu$ is employed~\cite{Ursem1999,Thomsen2004}. 
Another measure is the averaged Hausdorff distance (AHD)~\cite{Schutze2012}
\[
\operatorname{AHD}(P) = 
\max\left\{\left(\frac{1}{\nu}\sum_{i=1}^{\nu} d_{\mathrm{nn}}(\vec{x}^*_i, P)^p\right)^{1/p}, \left(\frac{1}{N}\sum_{i=1}^{N} d_{\mathrm{nn}}(\vec{x}_i, X^*)^p\right)^{1/p} \right\}\;,
\]
by using $X^*$ as a reference set.
The function $d_{\mathrm{nn}}(\vec{x}, X)$ denotes the Euclidean distance of a point $\vec{x}$ to its nearest neighbor in a set of points~$X$.

\begin{table}[t]
\caption[Different magnitudes for the number of function evaluations]{Different magnitudes for the number of function evaluations $N_f$.}
\centering
\small
\begin{tabularx}{\columnwidth}{lL}
\toprule
Magnitude & Application\\
\midrule
$n \cdot 10^1$ & Initial designs in model-based optimization \cite{Jones1998}\\
$n \cdot 10^2$ & Expensive optimization\\
$n \cdot 10^3$ & \\
$n \cdot 10^4$ & Budget of the CEC 2005 competition \cite{Suganthan2005}\\
$n \cdot 10^5$ & Budget of the CEC 2013 niching competition \cite{CEC2013niching}\\
$n \cdot 10^6$ & Budget of the black-box optimization benchmark (BBOB) \cite{wp200901_2010}\\
\bottomrule
\end{tabularx}
\label{tab:magnitudes_function_evals}
\end{table}

\begin{table}[t]
\caption{Test problems used in the experiment.} 
\centering
\small
\begin{tabularx}{\columnwidth}{lRRRR}
\toprule
Problem name & Dim. $n$ & \#Local optima & \#Global optima & Ref.\\
\midrule
Shekel5 & 4 & 5 & 1 & \cite{Dixon1978} \\
Shekel7 & 4 & 7 & 1 & \cite{Dixon1978} \\
Shekel10 & 4 & 10 & 1 & \cite{Dixon1978} \\
Hartman3 & 3 & 3 & 1 & \cite{Dixon1978} \\
Hartman6 & 6 & 2 & 1 & \cite{Dixon1978} \\
Goldstein-Price & 2 & 5 & 1 & \cite{Dixon1978} \\
Branin & 2 & 3 & 3 & \cite{Dixon1978} \\
Vincent & 2 & 36 & 36 & \cite{CEC2013niching} \\
Vincent & 3 & 216 & 216 & \cite{CEC2013niching} \\
Modified Rastrigin & 4 & 48 & 1 & \cite{DebSah12}\\
Modified Rastrigin & 8 & 48 & 1 & \cite{DebSah12}\\
Six-hump camelback & 2 & 6 & 2 & \cite{CEC2013niching} \\
\bottomrule
\end{tabularx}
\label{tab:test_problems}
\end{table}

\begin{figure*}[t]
    \centering
\includegraphics[width=\textwidth]{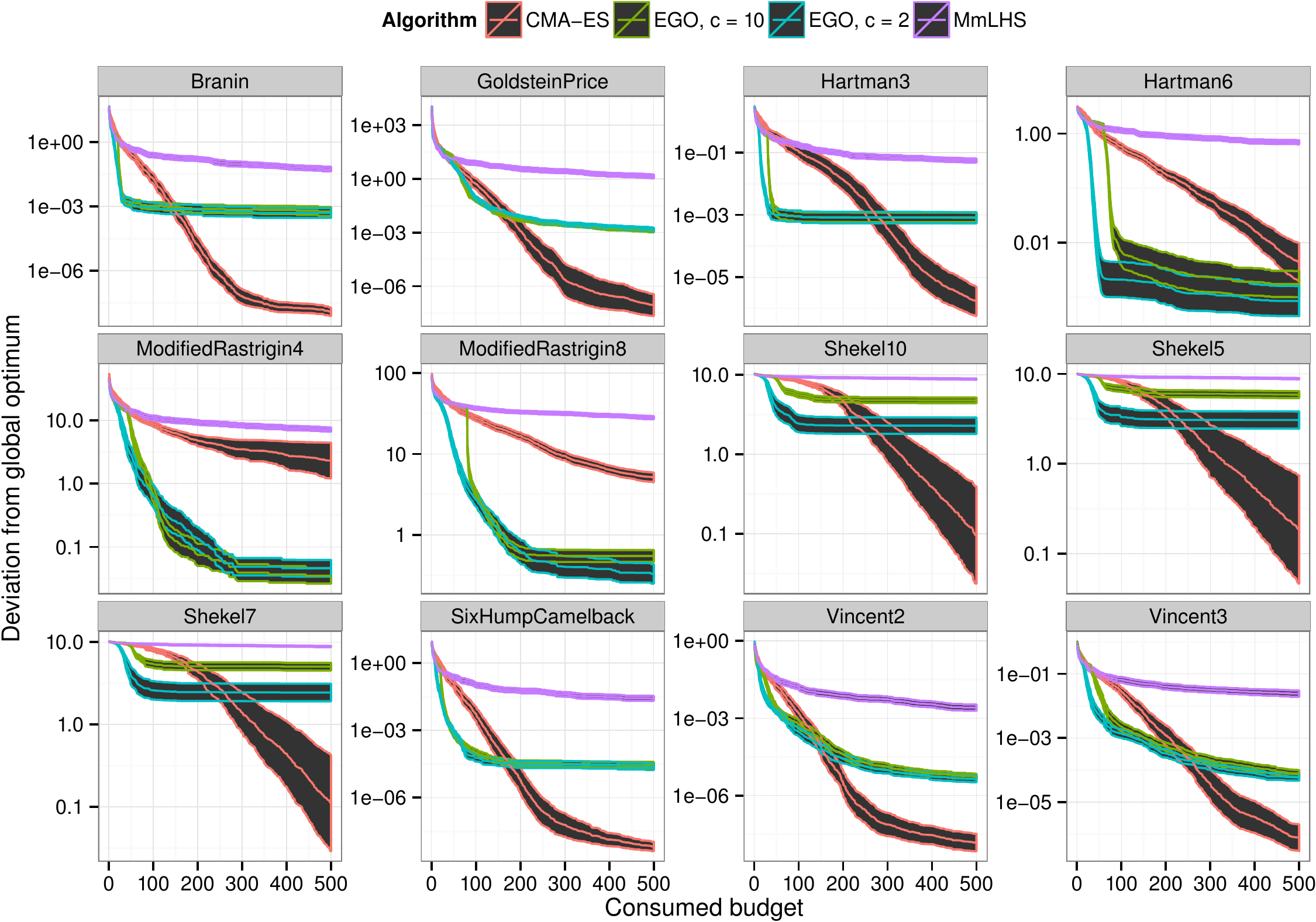}%
    \caption{Convergence graphs for approximating the global optimum. Note the logarithmic scale.}
	\label{fig:f_target_delta}
\end{figure*}

\section{Experiment}
\label{sec:experiment}

\noindent\textit{Research question:}
How does the assessment of optimization algorithms depend on the performance measurement in expensive optimization, i.\,e., are the results in global optimization different from multi-local optimization?

\noindent\textit{Pre-experimental planning:}
Table~\ref{tab:magnitudes_function_evals} shows how some budgets are associated with research areas and benchmarks. 
Measuring consumed resources simply as the number of objective function evaluations is generally deemed admissible if this number is small, because then the assumption of expensive function evaluations in relation to the overhead of an optimization algorithm holds.
For extremely large budgets, this is rather unlikely~\cite{Eiben2002}.
However, in expensive optimization, often very computationally demanding algorithms are used, so it is also an interesting question where the actual break-even point between two optimization algorithms is in terms of wall clock time.

\begin{figure*}[t]
    \centering
\includegraphics[width=\textwidth]{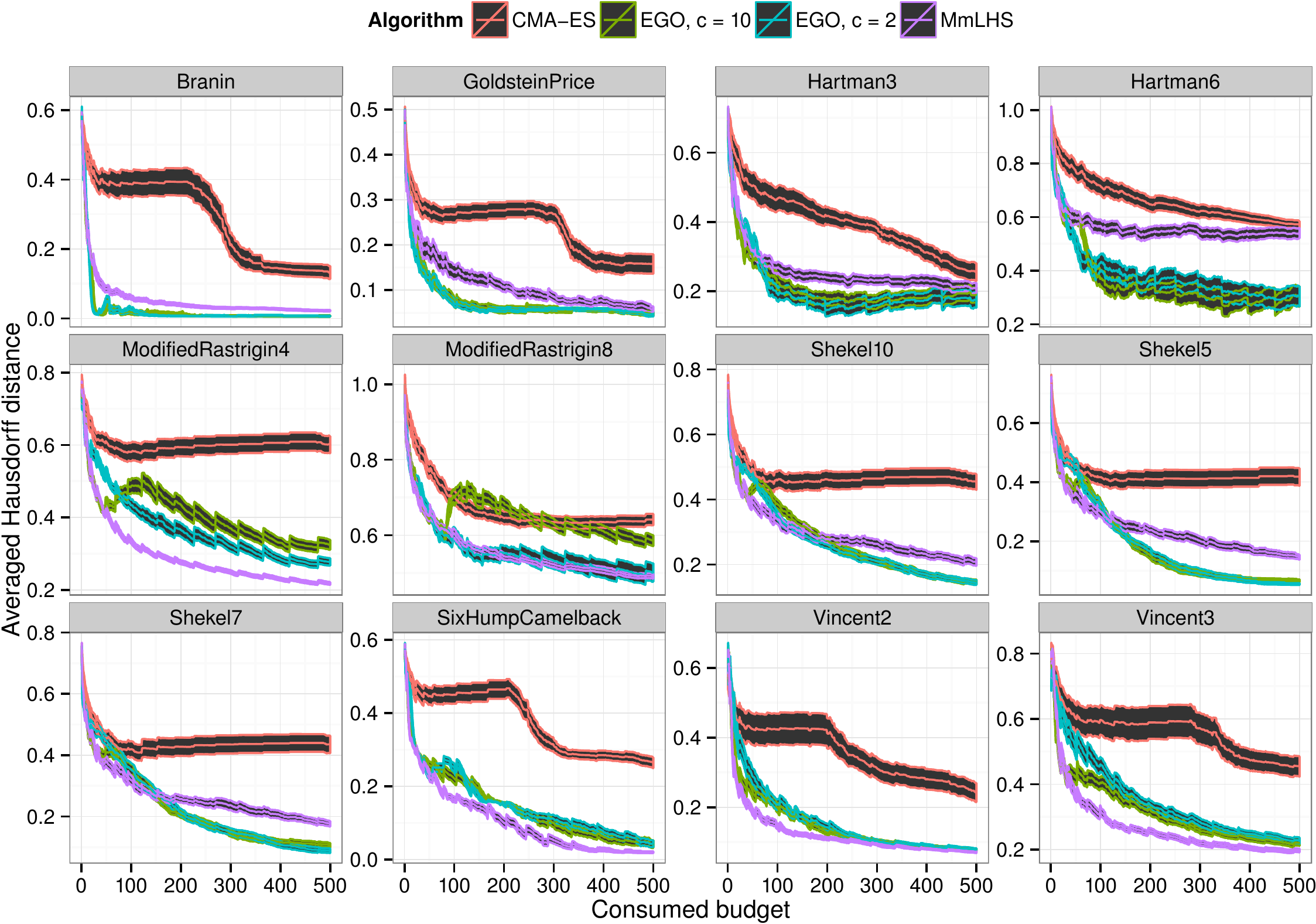}%
    \caption{Averaged Hausdorff distances (AHD) over the course of optimization.}
	\label{fig:ahd}
\end{figure*}

\noindent\textit{Task:}
We assume an anytime scenario for assessment, that is, the algorithms could be stopped at any time. 
We record three different performance measures over the course of optimization, namely the deviation from the global optimum, the peak ratio (PR), and the averaged Hausdorff distance (AHD). 
For PR, the position of an optimum is considered as approximated if a point is within a Euclidean distance of $0.01$ in the normalized search space (see setup below). 
AHD is used with an exponent of one.
For PR and AHD, the solutions up to the measuring point are filtered by topographical selection (TS), to stay close to a real-world scenario.
Topographical selection, originally proposed by~\cite{Toern1992}, contains a parameter $k$, specifying a number of neighbors.
To determine this parameter, we use the model
\begin{equation*}
k(n, N) = 0.215n + 0.74N^{1/2}\,,
\end{equation*}
depending on the dimension $n$ and the number of points $N$.
It was developed in~\cite{Wessing2016a} for random uniform samples.
Although the solution sets produced by the optimization algorithms are not uniformly distributed (except for MmLHS, see below), we feel certain that this is not a severe problem, as TS proved quite robust to changes in the distribution in previous experiments~\cite{Wessing2016a}.

\begin{figure*}[t]
    \centering
\includegraphics[width=\textwidth]{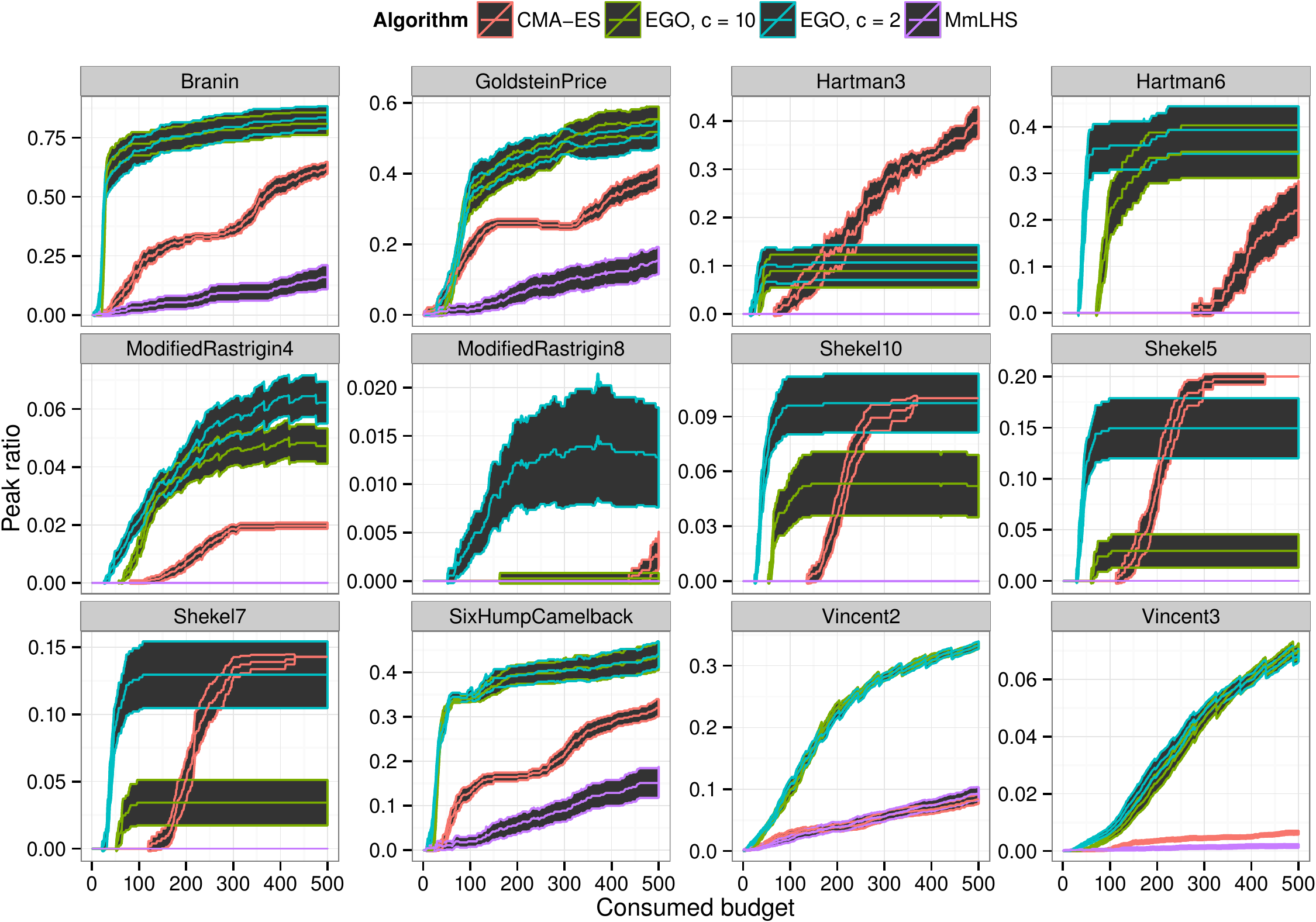}%
    \caption{Peak ratio (PR, $r = 0.01$) over the course of optimization.}
	\label{fig:pr_001}
\end{figure*}

\noindent\textit{Setup:}
Table~\ref{tab:test_problems} contains the test problems used in this experiment. 
They consist of the classic test set for global optimization by Dixon and Szeg\"o~\cite{Dixon1978}, and some problems taken from the 2013 niching competition~\cite{CEC2013niching}.
The former problems mostly contain fewer minima than the latter ones. 
However, the latter ones in part have other properties that make them easier, i.\,e., separability (modified Rastrigin) or no local optima (Vincent).
All problems have in common the rather low dimension, bound constraints, and the fact that positions of \emph{all} local and global minima are known. 
The last aspect is crucial for carrying out the assessment in the multi-local case with PR and AHD.
The search spaces are always normalized to the unit hypercube.

\begin{figure*}[t]
    \centering
\includegraphics[width=\textwidth]{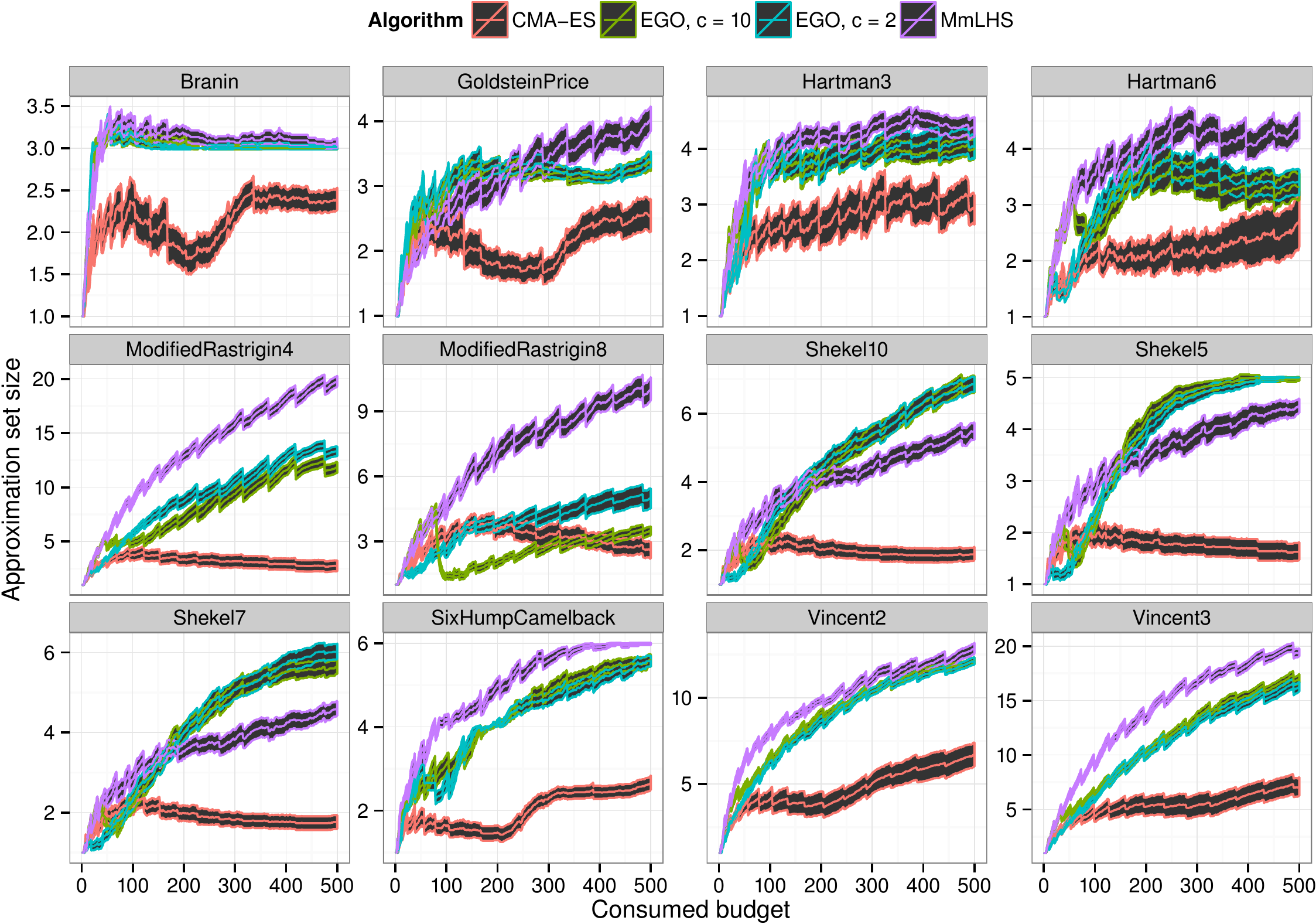}%
    \caption{The number of basins identified by topographical selection after each evaluation.}
	\label{fig:approx_set_size}
\end{figure*}

We compare three different algorithms, namely CMA-ES, EGO, and a maximin Latin hypercube sampling (MmLHS).
MmLHS acts as a representative of random search here.
Our MmLHS designs are produced on-the-fly by a greedy construction heuristic.
Thus, they are not exactly optimal according to the maximin-distance criterion, but possess a significantly higher uniformity than random uniform sampling (see~\cite[p.~58]{Wessing2015} for details).
For EGO, we try two variants, which only differ in the amount of function evaluations invested into the intial sample. 
The sample size is determined as $cn$, with $c = 2$ or $c = 10$, in accordance with Tab.~\ref{tab:magnitudes_function_evals}.
The initial sample for EGO is also drawn by MmLHS, thus MmLHS alone can be seen as a limiting case for EGO where the whole budget is spent on the initial sample.
EGO's Kriging model uses the power-exponential kernel
\begin{equation*}
\operatorname{corr}(\vec{x}_i, \vec{x}_j) = \exp\left(-\sum_{\ell=1}^n \theta_\ell |x_{i,\ell} - x_{j,\ell}|^{p} \right)\,.
\end{equation*}
For the kernel parameters, we require $-2 \leq \log_{10}(\theta_\ell) \leq 2$ and $0.5 \leq p \leq 2$.
As virtually all contemporary EGO implementations, our code differs from the algorithm in \cite{Jones1998} in the way the infill criterion is optimized.
Instead of a branch-and-bound approach, which consumes a lot of memory and restricts the kernel choice, we simply use CMA-ES, started once from the best of $100n$ uniformly distributed points.
Also the likelihood function for fitting the model is optimized with CMA-ES, based on recommendations in~\cite{Preuss2010}.

CMA-ES is the candidate in this test set with a strong focus on local search.
We use version 1.1.7 of the Python implementation\footnote{https://pypi.python.org/pypi/cma/}.
Its ``tolfun'' and ``tolfunhist'' stopping criteria are set to $10^{-3}$ and~$10^{-5}$, respectively, to stop really early and thus potentially have some budget left for starting another search. 
The starting points for CMA-ES are drawn by the maximin reconstruction algorithm, as recommended in~\cite{Wessing2015}. 
The initial step size is set to $0.15$.

In total, this experimental setup is chosen deliberately rather in favor of EGO than of CMA-ES, by including the test problems EGO was originally proposed for~\cite{Jones1998}. 
Additionally, CMA-ES is geared to being a very \emph{robust} black-box optimizer, so it is not necessarily the most efficient one on these low-dimensional, continuously differentiable problems~\cite{Rios2013}.

\begin{figure}[t!]
    \centering
\includegraphics[width=0.75\columnwidth]{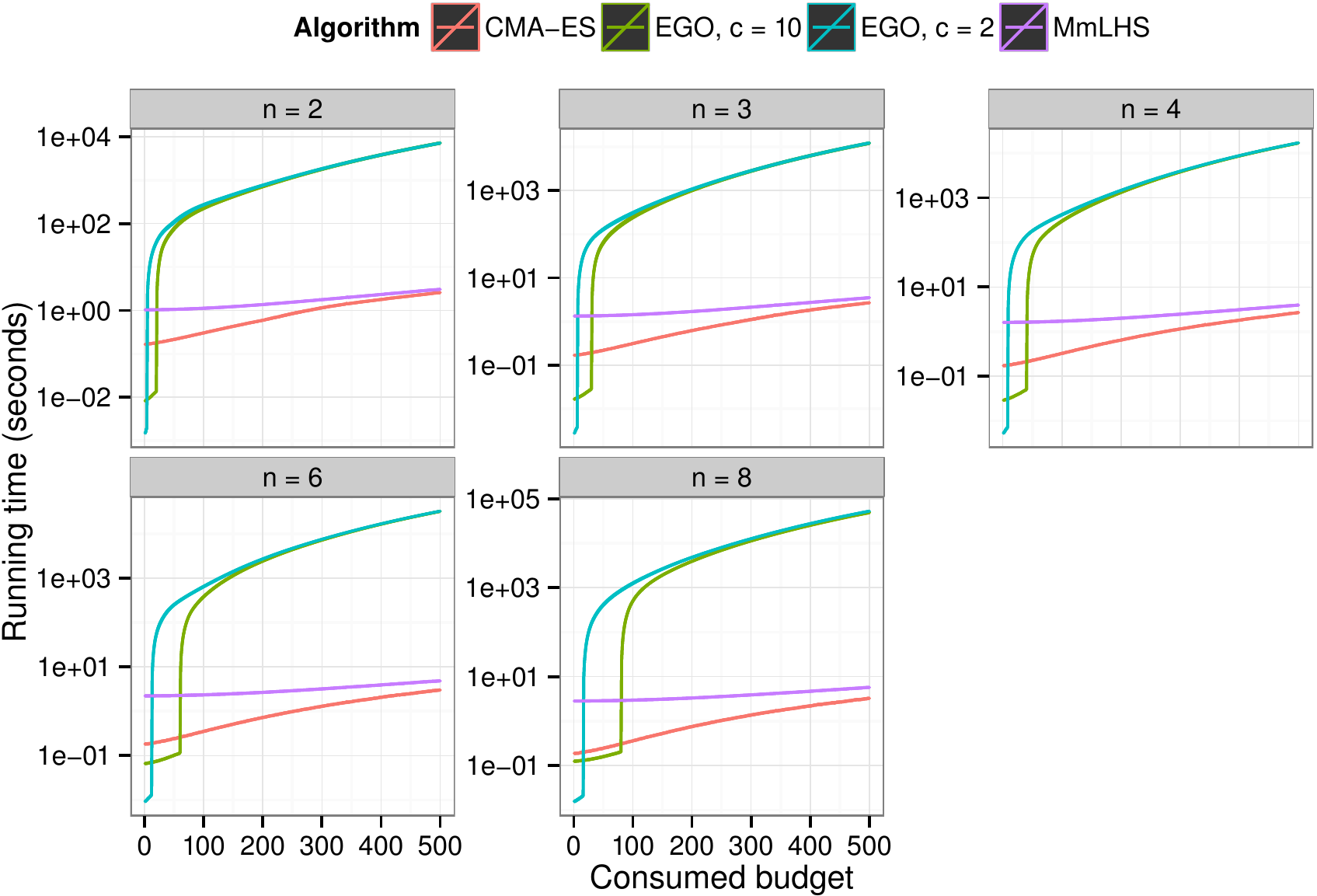}%
    \caption{The wall clock time required by the algorithms.}
	\label{fig:running_times}
\end{figure}

\noindent\textit{Results:}
Figures~\ref{fig:f_target_delta}, \ref{fig:ahd}, and~\ref{fig:pr_001} illustrate the development of the three indicators over the course of optimization.
Additionally, we report the number of selected solutions in Fig.~\ref{fig:approx_set_size}.
Figure~\ref{fig:running_times} shows the wall clock times of the algorithms. 
The curves contain the time for running the algorithm for the number of evaluations specified on the x-axis, plus the time for executing topographical selection once.
In all figures, the curves represent mean values over 75 stochastic replications, with a 95\% confidence interval for the mean under normality assumption.

\noindent\textit{Observations:}
Figure~\ref{fig:f_target_delta} shows that even under the quite restricted budget of 500 evaluations, CMA-ES significantly beats EGO on some problems EGO was developed for, if only the deviation from one global optimum counts. 
The variance is larger for CMA-ES, because some runs naturally only converge to local optima, but the average performance is clearly better.
However, EGO is always better than CMA-ES in terms of AHD and PR. 
On problems with a large number of optima, MmLHS obtains a still better AHD than EGO, but the optima are not approximated very well. 
Thus, EGO always has the better peak ratio.
With two exceptions, its peak ratio is also always better than that of CMA-ES.

Only few significant differences can be found between the choices $c = 2$ and $c = 10$ for EGO, but the results seem to be slightly in favor of $c = 2$.

Figure~\ref{fig:approx_set_size} shows that the number of selected solutions is mostly nicely correlated with the problems' actual number of optima. 
On Branin and Shekel5, even the correct number of optima is reliably identified towards the end of the optimization.
However, the approximation quality does not satisfy the PR criterion for all optima.

The running times of the algorithms in wall clock time naturally differ by orders of magnitude (see Fig.~\ref{fig:running_times}). 
The curve for MmLHS is slightly misleading, as the whole 500-point sample was computed in advance. 
So the curve begins with this cost and subsequently adds the cumulative cost of the test functions plus the cost for one topographical selection. 
In reality, MmLHS is of course always the cheapest algorithm in this setting, because one would not sample more points than one wants to evaluate.

\noindent\textit{Discussion:}
The results show that EGO represents an intermediate strategy between CMA-ES and MmLHS regarding the exploration-exploitation trade-off. 
It has the potential to detect several attraction basins with quite small budgets, but lacks the ability to approximate the corresponding optima with high precision. 
The $f_\Delta$ may stagnate for several hundred function evaluations. 
On the other hand, performance measures from multimodal optimization do often keep improving during this time.

The bends in the curves of CMA-ES in Fig.~\ref{fig:f_target_delta} are probably due to the aggressive stopping criteria, which prevent the algorithm from approximating the global optimum to a higher precision.
This is the only explanation on problems as Branin or Vincent, which only contain global optima, and where we would expect a linear convergence behavior otherwise. 
However, this is not to be seen as a drawback, as we deliberately chose this setting to obtain a better global search, and the restarts do clearly improve other measures (see Figs.~\ref{fig:ahd} and \ref{fig:pr_001}).
Tuning the initial step size might improve the CMA-ES performance  slightly more. 
However, note that any improvement of CMA-ES would further strengthen the support for our hypothesis, so the omission is not critical.

\section{Conclusions}
We showed that \emph{efficient global optimization} (EGO) is in fact not always the best algorithm for global optimization, i.\,e., the application it was originally designed for, except for extremely small budgets of approximately up to $200$ function evaluations. 
By our experimental setup, this statement is restricted to rather low-dimensional ($2 \leq n \leq 8$), smooth, and generally well-behaved objective functions. 
However, as higher-dimensional, more multimodal, and/or non-continuous functions would pose even more difficulties to the meta-modelling, and other optimization algorithms as, e.\,g., CMA-ES are inherently more robust to such difficulties, because they use less assumptions about the problem, the statement might be extended to broader settings in the future. 
Of course more sophisticated sequential model-based optimization algorithms do already exist, and may not share some of the basic EGO's weaknesses in global optimization, but our point is that EGO is actually fairly well-suited for the slightly different problem definition of finding multiple optima, when the optimization problem is expensive. 
Thus, combining it with a suitable basin identification heuristic makes it a strong competitor in this domain.

In future work, we shall look again at the basin identification mechanism and find better values or heuristics for the $k$ parameter, or even another algorithm altogether. 
Also, additional model-based optimization algorithms may be tested to find out if there exists an approach that is competitive in both the global and multi-local optimization case.
Finally, the methods shall be thoroughly benchmarked on a larger set of problems in order to make stronger claims on their strengths and weaknesses.

\end{document}